\documentclass[12pt,twoside,reqno]{amsart}

\usepackage{amsmath,amssymb,amsthm}
\usepackage{fancyhdr}
\usepackage{epsfig,color,colordvi} 

\usepackage[numbers,sort&compress,square,comma]{natbib}

\theoremstyle{definition}

\theoremstyle{remark}
\newtheorem{remark}{\sc Remark}

\theoremstyle{plain}
\newtheorem{theorem}{\sc Theorem}
\newtheorem{lemma}{\sc Lemma}
\newtheorem{proposition}{\sc Proposition}

\newcommand{\be}{\begin{equation}}
\newcommand{\ee}{\end{equation}}
\newcommand{\nn}{\nonumber}

\def\e{\varepsilon}
\def\ind{\mathbf{1}}  

\def\cF{\mathcal{F}}  
\def\cG{\mathcal{G}}

\def\Etil{{\widetilde{E}}}
\def\Ptil{{\widetilde{P}}}

\def\Sbar{{\bar{S}}}

\def\eps{\varepsilon}

\def\fR{{{\rm I}\!{\rm R}}}
\def\fN{{{\rm I}\!{\rm N}}}
\def\fZ{{{\rm Z}\mkern-5.5mu{\rm Z}}}

\def\fQ{{{\rm Q}\kern-.65em {}^{{}_/ }\,}}
\def\fQQ{ {{\rm Q}\kern-.57em \scriptscriptstyle{}^{]\kern.055em[}\,}}
\def\fP{{I\!\!P}}
\def\fE{{I\!\!E}}

\def\w{\omega}

\def\ord{\kern0.1em o\kern-0.02em{}_{\ds\breve{}}\kern0.1em}
\def\Ord{\kern0.1em O\kern-0.02em{\ds\breve{}}\kern0.1em}
\def\ds{\displaystyle}

\def\fmonth{\ifcase\month\or Jan\or Feb\or Mar\or Apr
\or May\or Jun\or Jul\or Aug\or Sep
\or Oct\or Nov\or Dec\fi\ }
\def\mmddyyyy{\the\month.\the\day.\the\year}
\def\ddmmyyyy{\the\day.\the\month.\the\year}
\def\Mddyyyy{\fmonth~\the\day,~\the\year}

\def\R{\fR}
\def\N{\fN}
\def\Z{\fZ}

\def\P{\fP}
\def\E{\fE}

\providecommand{\abs}[1]{\left\vert#1\right\vert}
\providecommand{\norm}[1]{\left\Vert#1\right\Vert}

\numberwithin{equation}{section}

\allowdisplaybreaks[1]

\begin{document}


\author{F.~Rassoul-Agha}   
\address{Mathematical Biosciences Institute, Ohio State University, 
Columbus, OH 43210}
\email{firas@math.ohio-state.edu}
\urladdr{www.math.ohio-state.edu/$\sim$firas}
\author{T.~Sepp\"al\"ainen}  
\address{Mathematics Department, University of Wisconsin-Madison, Madison, WI 53706}
\email{seppalai@math.wisc.edu}
\urladdr{www.math.wisc.edu/$\sim$seppalai}
\thanks{T.~Sepp\"al\"ainen was partially supported by
National Science Foundation grant DMS-0402231.}

\date{\today}





\begin{abstract}
We consider a discrete time random walk in a space-time i.i.d.~random
environment.
We use a martingale approach to show that
the walk is diffusive in almost every fixed environment.
We improve on existing results by proving an
invariance principle and considering environments with an annealed $L^2$ drift.
We also state an a.s.~invariance principle for random walks in
general random
environments whose hypothesis requires a subdiffusive bound on the
variance of the quenched mean,
under an ergodic invariant measure for the environment chain.
\end{abstract}




\title[An a.s.~invariance principle for RW in space-time RE]
{An almost sure invariance principle for random walks in a space-time random environment}

\maketitle


\def\X{{\hat X}}

\section{Introduction}
\label{intro} Random walk in a random environment is one of the
basic models of the field of disordered systems of particles. In
this model, an environment is a collection of transition
probabilities $\w=(\pi_{xy})_{x,y\in\Z^d}\in{\mathcal P}^{\Z^d}$
where
 ${\mathcal P}=\{(p_z)_{z\in\Z^d}\in[0,1]^{\Z^d}:\sum_z p_z=1\}$.
Let us denote by $\Omega={\mathcal P}^{\Z^d}$ the space of all
such transition probabilities. The space $\Omega$ is equipped with
the canonical product $\sigma$-field $\mathfrak{S}$ and with the
natural shift $\pi_{xy}(T_z\w)=\pi_{x+z,y+z}(\w)$, for $z\in\Z^d$.
On the space of environments $(\Omega,\mathfrak{S})$, we are given
a certain $T$-invariant probability measure $\P$ with
$(\Omega,\mathfrak{S},(T_z)_{z\in\Z^d},\P)$ ergodic. We will say
that the environment is i.i.d.~when $\P$ is a product measure in
the sense that the random probability vectors
$(\pi_{x,y})_{y\in\Z^d}$ are i.i.d.~over distinct sites $x$. 
We will denote by $\E$ the expectation under $\P$.
Let us now describe the process.

First, the environment $\w$ is chosen
from the distribution $\P$. Once this is done it remains fixed for
all times. The random walk in environment $\w$,
starting at $z$, is then the canonical
Markov chain $\X=(X_n)_{n\geq0}$ with state space $\Z^d$ and
satisfying
\begin{align*}
P_z^\w(X_0=z)&=1,\hfill\cr
P_z^\w(X_{n+1}=y|X_n=x)&=\pi_{xy}(\w).\hfill
\end{align*}
The process $P_z^\w$ is called the {\sl quenched law}.
The {\sl joint annealed law} is then
$$P_z(d\X,d\w)=P_z^\w(d\X)\P(d\w),$$
and $P_z(d\hat X,\Omega)$ is
the {\sl marginal annealed law} or, by {\it abus de langage},
just the {\sl annealed law}. 
We will use $E_z$ and
$E_z^\w$ for the expectations under, respectively, $P_z$
and $P_z^\w$. 

We begin by considering a special type of random environment.
Namely, assume that $d=\nu+1\geq2$ and let $\{e_i\}_{i=1}^d$ be
the canonical basis of $\R^d$. Assume then that $\P$ is
i.i.d.~and, $\P$-a.s.,
\begin{align}
\label{spacetime}
\pi_{0,e_1+z}=0,\hbox{ if }z\not\in
E=\{x\in\Z^d:x\cdot e_1=0\}\sim\Z^{\nu}.
\end{align}
Condition (\ref{spacetime}) says that at each step the first
coordinate increases by one deterministically. One of the reasons
for considering such a model comes from the fact that if one views
it as a random walk in a $\nu$-dimensional space-time
i.i.d.~random environment, with $\R e_1$ being the time axis, then
it turns out to be a dual process to some surface growth
processes. See, for example, \cite{FF}.

Clearly, the annealed process, in this case, is equivalent to
$(ne_1+Y_n)_{n\geq0}$, where
$Y_n$ is a homogeneous Markovian random walk on $E$, with
transitions
$$(p(x,x+z)=p(0,z)=\E(\pi_{0,e_1+z}))_{x,z\in E}.$$
If $p$ has a first moment, the annealed walk on $\Z^d$ has a law of large
numbers with velocity $$v=e_1+\sum_{z\in E}p(0,z)z.$$
If, furthermore, it has a second moment, one then has an annealed
invariance principle with diffusion matrix
$${\mathfrak D}=\sum_{z\in E}(e_1+z-v)(e_1+z-v)^t p(0,z).$$
Since any quenched central limit type result would imply an
annealed one, one can then see that a second moment condition on
$p$ has to be assumed, if one wants to prove an a.s.~invariance
principle for this model. On the other hand, if
$\P(\sup_z\pi_{0z}=1)=1$, then the quenched walk becomes
deterministic and a central limit is out of question. This
justifies our hypothesis (M for moment, E for ellipticity):

\newtheorem*{HypME}{Hypothesis (ME)}{\bf}{\it}
\begin{HypME}
The measure $\P$ satisfies
the following condition
\begin{align*}
\sum_{z\in E}|z|^2\E(\pi_{0,e_1+z})<\infty
\hbox{ and }\P(\sup_{z\in E}\pi_{0,e_1+z}<1)>0.
\end{align*}
\end{HypME}

Define now, for $t\ge0$ and a given $(X_n)_{n\geq0}$,
\begin{align}
\label{Bndef}
B_n(t)=\frac{X_{[nt]}-[nt]v}{\sqrt{n}}\hbox{ and }
\widetilde B_n(t)=\frac{X_{[nt]}-E_0^\w(X_{[nt]})}{\sqrt{n}}.
\end{align}
Here, for $x\in\R$, $[x]=\max\{n\in\Z: n\le x\}$.
For a closed interval $I\subset[0,\infty)$ denote by $D_{\R^d}(I)$
the space of right continuous functions on $I$, taking values in $\R^d$,
and having left limits.
The space $D_{\R^d}(I)$ is endowed with the usual
Skorohod topology \cite{EK}.
For $\w\in\Omega$,
let $Q_n^\w$, respectively $\widetilde Q_n^\w$, denote the distribution of $B_n$,
respectively $\widetilde B_n$, induced by $P_0^\w$, on the
Borel sets of $D_{\R^d}([0,\infty))$.
In Section \ref{RWSTRE} we will prove the following theorem:

\begin{theorem}
\label{cltdirected} Let $d\geq2$ and consider a random walk in an
i.i.d.~random environment satisfying {\rm (\ref{spacetime})} and
Hypothesis {\rm (ME)}. Then, for $\P$-a.e. $\w$, the distribution
$Q_n^\w$ converges weakly to the distribution of a Brownian motion
with diffusion matrix given by $\mathfrak D$. Moreover,
$n^{-1/2}\max_{k\leq n}|E_0^\w(X_n)-nv|$ converges to $0$,
$\P$-a.s. and, therefore, the same invariance principle holds also
for $\widetilde Q_n^\w$.
\end{theorem}

\begin{remark}
\label{BMdef}
For a symmetric, non-negative definite $d\times d$ matrix $\Gamma$,
a Brownian motion with diffusion matrix $\Gamma$ is the $\R^d$-valued
process $(W(t))_{t\ge0}$
such that $W(0)=0$, $W$ has continuous paths, independent increments,
and for $s<t$ the $d$-vector $W(t)-W(s)$ has Gaussian distribution
with mean zero and covariance matrix
$(t-s)\Gamma$. If the rank of $\Gamma$ is $m$,
one can produce such a process by
finding a $d\times m$ matrix $\Lambda$ such that
$\Gamma=\Lambda\Lambda^t$, and by defining $W(t)=\Lambda B(t)$
where $B$ is an  $m$-dimensional standard Brownian motion.
\end{remark}

In space-time product random environments the invariance principle under the
annealed $P_0$ is just Donsker's classical invariance principle. But in
general random environments even the annealed invariance
principle is far from immediate. For recent
results in this direction see \cite{czclt} and \cite{effective} and the
references therein.
Moreover, the switch from annealed central limit type results to
quenched ones for random walks in random environments
is a hard problem that has been subject to a very slow progress.

Quenched results have been proved under specialized assumptions.
For example, if the random walk is balanced (i.e.~$\pi_{0x}=\pi_{0,-x}$),
then $X_n$ becomes a martingale under the quenched measure, and one can
show a quenched invariance principle \cite{lawler}.
On the other hand,
if the random walk has a sufficiently high-dimensional
simple symmetric random walk part,
then one can use a natural regeneration structure that arises
\cite{cutpoints}. If the random
environment is a small perturbation of the simple symmetric random
walk, then a quenched invariance principle has been proved using
renormalization techniques; see \cite{BK} and the recent \cite{cltrwrenoise}.
Even in the basic reversible case of random walks among random conductances,
quenched invariance principles have only recently been shown to hold
\cite{conductance}. There the authors
adopt the approach developed by Kipnis and Varadhan \cite{KV} for reversible
Markov chains to get a central limit theorem; see Section \ref{corrector}
below. Also, see \cite{olla} and \cite{noteclt} for a review on the
approach in \cite{KV}.
For a non-reversible setting where \cite{KV} has also been useful,
see \cite{persistentrwre}.
Notably, as it is usually the case in this field,
the invariance principle does not immediately follow from the fixed time
central limit theorem. In \cite{conductance}, recent Gaussian estimates are used
to perform the transition.

In the case of a space-time product environment, the central limit
version for $B_n(1)$ has been known in the case of small noise;
see \cite{qclt-noise} and \cite{stannat}. In the recent
\cite{qclt-nonoise} the authors remove the ``small noise''
assumption and prove the fixed time central limit theorem under
the hypothesis that $\sup_\w\pi_{0z}(\w)$ has an exponential
moment, as opposed to our second moment assumption in Hypothesis
(ME). If one assumes on top of Hypothesis (ME) that $\P$-a.s.,
$\pi_{0,e_1+z}=0$, if $z\not\in\{e_i,-e_i\}_{i\geq2}$, and
positive otherwise, then the invariance principle for the process
$B_n$ in dimension $d\geq4$ has been proved in \cite{dynstat}.

Our arguments are different from \cite{qclt-nonoise}. The ideas
in \cite{dynstat} and \cite{conductance} are somewhat related in spirit.
Our approach is based on adaptations of the well-known
Kipnis-Varadhan method \cite{KV} to non-reversible situations
developed by Maxwell and Woodroofe \cite{MW} and Derriennic
and Lin \cite{DL}. Maxwell and Woodroofe use fairly concrete
probabilistic reasoning. We will refer directly to their paper for
some preliminary steps in our proof. The approach of Derriennic
and Lin is more abstract and powerful, cast in the framework of
Banach space contractions, and ultimately produces stronger
results. We apply their results to conclude our proof. All this
will be described in Section \ref{martingales} where we prove a
result that holds for general random environments.

In the course of this paper $\w$, $\w_0$, and $\w_1$ will denote
generic elements of $\Omega$. 
We will write $A^t$ for the transpose of a vector or
matrix $A$. An element of $\R^d$ is regarded as a $d\times 1$
matrix, or column vector. The set of whole numbers
$\{0,1,2,\cdots\}$ will be denoted by $\N$.

\section{Quenched invariance principle for general random environments
under moment hypotheses}
\label{martingales}
In this section we consider the general random walk in a random environment
 as formulated
in the first two paragraphs of the Introduction. In particular,
the special structure of assumption (\ref{spacetime}) and
Hypothesis (ME) are not assumed in this section.

First, let us define the drift $$D(\w)=E_0^\w(X_1)=\sum_z z\pi_{0z}(\w).$$
For a bounded measurable function $h$ on $\Omega$, define
$$\Pi h(\w)=\sum_z\pi_{0z}(\w)h(T_z\w).$$
In fact, $\Pi-I$ defines the generator of the Markov process of the environment,
as seen from the particle. This is the process on $\Omega$ with transitions
$$\pi(\w,A)=P_0^\w(T_{X_1}\w\in A).$$
In this section, we will assume there exists a probability measure
$\P_\infty$ on $\Omega$ that is invariant for the transition
$\Pi$. Then, the operator $\Pi$ can be extended to a contraction
on $L^p(\P_\infty)$, for every $p\in[1,\infty]$. We will use the
notation $\E_\infty$ for the corresponding expectation. When the
initial distribution is $\P_\infty$, we will denote this Markov
process by $\widetilde P_0^\infty$. We will also write
$P_0^\infty(d\X,d\w)=P_0^\w(d\X)\P_\infty(d\w)$ and $E_0^\infty$
for the expectation under $P_0^\infty$. Note that
 $\widetilde P_0^\infty$ is the probability measure induced by $P_0^\infty$ and
$(T_{X_n}\w)$ onto $\Omega^\N$. With this notation, the measure
$$\mu_2^\infty(d\w_0,d\w_1)=\pi(\w_0,d\w_1)\P_\infty(d\w_0)$$
describes the law of $(\w,T_{X_1}\w)$, under $P_0^\infty$.

We will assume in this section that
$D\in L^2(\P_\infty)$.
Next, for $\eps>0$, let $h_\eps$ be the solution of
$$(1+\eps)h_\eps-\Pi h_\eps=g=D-v,$$
where $v=\E_\infty(D)$. In fact, one can write:
$$h_\eps=\sum_{k=1}^\infty(1+\eps)^{-k}\Pi^{k-1}g\in L^2(\P_\infty).$$
Define
$$H_\eps(\w_0,\w_1)=h_\eps(\w_1)-\Pi h_\eps(\w_0).$$
Then one has the following theorem.

\begin{theorem}
\label{MW}
Let $d\geq1$ and
let $\P_\infty$ be any probability measure on $(\Omega,{\mathfrak S})$
that is invariant under $\Pi$ and ergodic for the Markov process
on $\Omega$ with generator $\Pi-I$. Assume that
$\sum_z|z|^2\E_\infty(\pi_{0z})<\infty$.
Assume also that there exists an $\alpha<1/2$ such that
\begin{align}
\label{cond}
\sqrt{\E_\infty\left(\abs{E_0^\w(X_n)-nv}^2\right)}=
\norm{\sum_{k=0}^{n-1}\Pi^k g}_2=\Ord(n^\alpha),
\end{align}
where $\norm{\cdot}_p$ is the $L^p(\P_\infty)$-norm.
Then we get the following conclusions: The limit
\begin{align}
\label{H}
H=\lim_{\eps\rightarrow0^+}H_\eps
\end{align}
exists in $L^2(\mu_2^\infty)$.
For $\P_\infty$-a.e.~$\w$, the distribution $Q_n^\w$ of the process $B_n$ defined
in (\ref{Bndef}) converges weakly to the
distribution of a Brownian motion with diffusion matrix
\begin{align}
\label{matrix}
E_0^\infty\bigl[(X_1-D(\w)+H(\w,T_{X_1}\w))(X_1-D(\w)+H(\w,T_{X_1}\w))^t\bigr],
\end{align}
as defined in Remark \ref{BMdef}. Moreover,
 $n^{-1/2}\max_{k\leq n}|E_0^\w(X_k)-kv|$ converges to $0$,
$\P_\infty$-a.s. and, therefore, the same invariance principle
holds for $\widetilde Q_n^\w$.
\end{theorem}

\begin{remark}
\label{1dcase}
When (\ref{cond}) only holds for $\alpha=1/2$, e.g. when $d=1$,
the quenched central limit
theorem may only hold with random centering $E_0^\w(X_n)$ and not with
deterministic centering $nv$. See Examples 3 and 4 and Proposition 1
of \cite{forbidden}.
\end{remark}

\begin{proof}
The proof uses essentially the strategy of \cite{DL} which
is an extension of some of the results of \cite{MW}, to which in
turn we will refer the reader for part of the calculations.

First let us give a few more definitions. For $\eps>0$, let
$$M_n^\eps=\sum_{k=0}^{n-1}H_\eps(T_{X_k}\w,T_{X_{k+1}}\w),~
\bar X_n=X_n-\sum_{k=0}^{n-1}D(T_{X_k}\w).$$
Furthermore, define $$S_n^\eps=\sum_{k=0}^{n-1}h_\eps(T_{X_k}\w),
~R_n^{\eps}=h_\eps(\w)-h_\eps(T_{X_n}\w),$$
so that $$X_n-nv=\bar X_n+M_n^\eps+\eps S_n^\eps+R_n^\eps.$$
Now, we proceed with the proof.
The existence of the limit in (\ref{H}) follows from Proposition 1 of
\cite{MW}. Thus, if one defines
$$M_n=\sum_{k=0}^{n-1}H(T_{X_k}\w,T_{X_{k+1}}\w),$$
then, for $\P_\infty$-almost every $\w$,  $(M_n)_{n\geq1}$ is a $P_0^\w$-square
integrable martingale relative to the filtration
$\{{\mathcal F}_n=\sigma(X_0,\cdots,X_n)\}_{n\geq0}$.
It also follows from Lemma 1 of \cite{MW} that
one has $\norm{h_\eps}_2=\Ord(\eps^{-\alpha})$. Define the error by
\begin{align}
\label{error}
R_n=X_n-nv-\bar X_n-M_n=M^\eps_n-M_n+\eps S^\eps_n+R^\eps_n.
\end{align}
Corollary 4 of \cite{MW} shows that
\begin{align}
\label{cond1}
E_0^\infty(\abs{R_n}^2)=\Ord(n^{2\alpha}).
\end{align}
%
%

Let $M_n^*(t)=n^{-1/2}(\bar X_{[nt]}+M_{[nt]})$. $(M_n^*(t))_{0\le
t\le1}$ converges weakly, under $P_0^\w$ for
$\P_\infty$-a.e.~$\w$, to a Brownian motion with diffusion matrix
as in (\ref{matrix}). This follows from a vector-valued version of
a well-known invariance principle for martingales. For the convenience
of the reader we provide a proof of it in the appendix.
The limits needed as hypotheses for this invariance principle
follow from ergodicity and the square-integrability of $M_1$ and
$X_1$. In turn, the assumption $\sum\abs{z}^2\E_\infty(\pi_{0z})<\infty$ 
guarantees that $X_1$ and $M_1$ are square-integrable for 
$\P_\infty$-a.e. $\w$. 

We have
$$
\sup_{0\le t\le 1} \abs{B_n(t)-M_n^*(t)} \le
n^{-1/2}\max_{k\le n}\abs{R_k}.
$$
Therefore, the invariance principle for $(B_n(t))_{0\le t\le1}$
will follow once we show that
\begin{align}
\label{R}
n^{-1/2}\max_{k\le n}\abs{R_k}
\mathop{\longrightarrow}_{n\rightarrow\infty}0,~
\hbox{in }P_0^\w\hbox{-probability, for }\P_\infty\hbox{-a.e. }\w.
\end{align}
By (\ref{error}), $E^\w_0(R_n)=E^\w_0(X_n)-nv$. Hence the
invariance principle for the process $(\widetilde B_n(t))_{0\le
t\le1}$ will follow once we show that
\begin{align}
\label{ER}
n^{-1/2}\max_{k\le n}\abs{E_0^\w(R_k)}
\mathop{\longrightarrow}_{n\rightarrow\infty}0,
\hbox{ for }\P_\infty\hbox{-a.e. }\w.
\end{align}

To prove (\ref{R}) and (\ref{ER}) we apply the theory of
``fractional coboundaries'' of Derriennic and Lin \cite{DL2}. The
first application is to the shift map $\theta$ on the sequence
space $\Omega^{\N}$ which is a contraction on the space
$L^2(\widetilde P^\infty_0)$.  On $\Omega\times\Omega$ define
first $f(\w_0,\w_1)=g(\w_0)-H(\w_0,\w_1)$. Then,
$P_0^\infty$-a.s.
$$R_n = \sum_{k=0}^{n-1}f(T_{X_k}\w,T_{X_{k+1}}\w).$$
 For sequences $\bar\w=(\w^{(i)})_{i\in\N}\in\Omega^{\N}$ define
$F(\bar\w)=f(\w^{(0)},\w^{(1)})$ and
 $$\widetilde R_n=\sum_{k=0}^{n-1} F\circ\theta^k.$$ Then
$F\in L^2(\widetilde P^\infty_0)$ and the process $(\widetilde R_n)_{n\ge1}$ has the
same distribution under $\widetilde P_0^\infty$ as the process $(R_n)_{n\ge1}$ has under
$P_0^\infty$.

Condition (\ref{cond1}) shows that the assumptions of Theorem 2.17
of \cite{DL2} are satisfied. The conclusion is that
$F\in(I-\theta)^\eta L^2(\widetilde P_0^\infty)$, for any
$\eta\in(0,1-\alpha)$. Since $\alpha<1/2$, we can find such an
$\eta\in(1/2,1-\alpha)$. But then (i) in Theorem 3.2 of \cite{DL2}
implies that $n^{-1/2}\widetilde R_n$ converges to $0$,
$\widetilde P_0^\infty$-a.s. This implies that $n^{-1/2}R_n$
converges to $0$, $P_0^\infty$-a.s. In other words, $n^{-1/2}R_n$
converges to $0$, $P_0^\w$-a.s., for $\P_\infty$-a.e.~$\w$. From
this, (\ref{R}) is immediate.

To prove (\ref{ER}) apply the same results from \cite{DL2} to the
contraction $\Pi$ on $L^2(\P_\infty)$.  Because
$$E_0^\w(R_n)=E_0^\w(X_n)-nv=\sum_{k=0}^{n-1}\Pi^k g$$
repeating the above argument with $\theta$ replaced by $\Pi$ and
(\ref{cond1}) replaced by (\ref{cond}) proves (\ref{ER}).

Once  we have the invariance principle on $D_{\R^d}([0,1])$, the
identities
$$B_n(t)=\sqrt{N} B_{nN}(t/N)\hbox{ and }
\widetilde B_n(t)=\sqrt{N}\widetilde B_{nN}(t/N)$$
show that on $D_{\R^d}([0,N])$ $B_n$ and $\widetilde B_n$ converge
to the process $\sqrt{N} W(t/N)$ which  is the same as $W$.
Then weak convergence on each $D_{\R^d}([0,N])$ implies
weak convergence on $D_{\R^d}([0,\infty))$.
The proof of Theorem \ref{MW} is complete.
\end{proof}

\section{On the corrector function}
\label{corrector}
In this section, we will prove an interesting property of the error
$R_n$ defined in (\ref{error}).  This property will show
the key difference between the one and multi-dimensional cases.
Although we will not make use of this property, it establishes a connection
with other existing ways of approaching the problem; see \cite{conductance}, for
instance.

First, define the functions
\begin{align*}
f_\eps(\w_0,\w_1)&=g(\w_0)-H_\eps(\w_0,\w_1)-\eps h_\eps(\w_0),\\
f(\w_0,\w_1)&=g(\w_0)-H(\w_0,\w_1).
\end{align*}
A calculation shows us that
\begin{align*}
f_\eps(\w_0,\w_1)&=h_\eps(\w_0)-h_\eps(\w_1),\\
R_n^\eps&=\sum_{k=0}^{n-1}f_\eps(T_{X_k}\w,T_{X_{k+1}}\w),\\
R_n&=\sum_{k=0}^{n-1}f(T_{X_k}\w,T_{X_{k+1}}\w).
\end{align*}
The following proposition establishes a ``co-cycle'' property of $f$.

\begin{proposition}
For $m,n\in\N$, let $\left((x_i)_{i=0}^n,(\tilde x_j)_{j=0}^m\right)$ be
``an admissible bridge'', i.e.~such that
$x_n=\tilde x_m=x$, and
$$\int P_0^\w(X_i=x_i,0\leq i\leq n)P_0^\w(X_j=\tilde x_j,0\leq j\leq m)
\P_\infty(d\w)>0.$$
Then, one has
\begin{align}
\label{claim}
\sum_{i=0}^{n-1}f(T_{x_i}\w,T_{x_{i+1}}\w)=
\sum_{j=0}^{m-1}f(T_{\tilde x_j}\w,T_{\tilde x_{j+1}}\w),~\P_\infty\hbox{-a.s}.
\end{align}
\end{proposition}

\begin{proof}
Fix $m,n$ and $x$, and let $Q_{0,x,n,m}^\w$ be the measure on
the space of double paths $((x_i)_{i=0}^n,(\tilde x_j)_{j=0}^m)$ defining two
independent random walks driven by the same environment $\w$, both starting at $0$
and ending at $x$ after, respectively, $n$ and $m$ steps. That is,
$$Q_{0,x,n,m}^\w(((x_i)_{i=0}^n,(\tilde x_j)_{j=0}^m))=
\prod_{i=0}^{n-1}\pi_{x_i x_{i+1}}(\w)
\prod_{j=0}^{m-1}\pi_{\tilde x_j\tilde x_{j+1}}(\w),$$
if $x_n=\tilde x_m=x,~x_0=\tilde x_0=0$, and $0$ otherwise.
Define also $$Q_{0,x,n,m}=\int Q_{0,x,n,m}^\w \P_\infty(d\w).$$ Then
one has, for each $\eps>0$,
$$\sum_{i=0}^{n-1}f_\eps(T_{X_i}\w,T_{X_{i+1}}\w)=
\sum_{j=0}^{m-1}f_\eps(T_{\tilde X_j}\w,T_{\tilde X_{j+1}}\w),~Q_{0,x,n,m}\hbox{-a.s.}$$
By (\ref{H}) each term above converges, in
$L^2(P_0^\infty)$, to the corresponding term in (\ref{claim}).
Note now that, although $Q_{0,x,n,m}$ is not a probability measure, one still
has $\frac{dQ_{0,x,n,m}}{dP_0^\infty}\leq1$.
Therefore, the convergence also happens in $L^2(Q_{0,x,n,m})$ and as a result
$$\sum_{i=0}^{n-1}f(T_{X_i}\w,T_{X_{i+1}}\w)=
\sum_{j=0}^{m-1}f(T_{\tilde X_j}\w,T_{\tilde X_{j+1}}\w),~Q_{0,x,n,m}^\w\hbox{-a.s.,
for }\P_\infty\hbox{-a.e. }\w.$$
Therefore, (\ref{claim}) holds for any admissible bridge.
\end{proof}

Let us, for simplicity, assume that all points in $\Z^d$ can be reached from
$0$ by some admissible path.
The above Lemma tells us then that if one defines the so-called
``corrector function'' as
$$\chi(x,\w)=\sum_{i=0}^{m-1}f(T_{x_i}\w,T_{x_{i+1}}\w),$$
where $(x_0,\cdots,x_m)$ is any admissible path from $0$ to $x$, then
\begin{align}
\label{cocycle}
R_n=\chi(X_n,\w)=\chi([nv],\w)
+\chi(X_n-[nv],T_{[nv]}\w).
\end{align}
Here, for $x\in\R^d$, $[\,\cdot\,]$
acts on each coordinate separately, to give $[x]$.
Relation (\ref{cocycle}) has a quite interesting implication.
The term $\chi([nv],\w)$ represents the fluctuations coming from the environment
itself.
When $d=1$, $\chi([nv],\w)$ is of order $\sqrt n$ and for a quenched invariance
principle one needs to consider $X_n-nv-\chi([nv],\w)$, i.e. to have a random
centering. See \cite{stflour} and Example 4 of \cite{forbidden} for more details.
However, when $d\geq2$ the walker can see ``more'' environments.
The quantity $\frac{\chi([nv],\w)}{\sqrt n}$ then vanishes at the limit and condition
(\ref{R}) has a chance to hold. See \cite{conductance} for a
result where control of the corrector is key to a quenched invariance principle
in a reversible setting.

\section{Space-time i.i.d.~random environments}
\label{RWSTRE}
In this section, we will prove Theorem \ref{cltdirected}
via an application of Theorem \ref{MW}.

\begin{proof}[of Theorem {\rm\ref{cltdirected}}]
For $n\geq0$, let
$$f_n(\w)=\sum_{x:x\cdot e_1=-n}P_x^\w(X_n=0).$$
By translation invariance of $\P$ one has $\E(f_n)=1$.
One then can check that $f_n$ is a martingale relative to the filtration
$$\{{\mathfrak S}_{-n}=\sigma((\pi_{xy})_y,~x\cdot e_1\geq-n)\}_{n\geq0}.$$
Therefore, there is a probability measure $\P_\infty$ such that
$$\frac{d\P_\infty{}_{|_{{\mathfrak S}_{-n}}}}{d\P_{|_{{\mathfrak S}_{-n}}}}=f_n.$$
An induction on $f_n$ shows also that $\P_\infty$ is invariant for $\Pi$.
%
Since $f_0=1$, we have that $\P=\P_\infty$ on ${\mathfrak S}_0$.
This will be of great use to us. On the one hand, since
$(\pi_{0z})_z$ is ${\mathfrak S}_0$-measurable, Hypothesis (ME)
implies that $D\in L^2(\P_\infty)$. On the other hand, on
$\mathfrak{S}_0$ the i.i.d.~structure carries over to $\P_\infty$.
Using this, we will show in the following lemma that $\P_\infty$
is also ergodic for the Markov process with generator $\Pi-I$.

\begin{lemma}
\label{ergodicity}
If $\P$ is i.i.d.~and satisfies {\rm (\ref{spacetime})}, then the invariant
measure $\P_\infty$, constructed above, is ergodic  for the Markov process
with generator $\Pi-I$.
\end{lemma}

\begin{proof}
Consider a bounded local function $\Psi$ on $\Omega$ that is measurable
with respect to $\sigma((\pi_{xy})_y,|x\cdot e_1|\leq K)$, for some integer
$K\geq0$.
Due to (\ref{spacetime}), $(\Psi(T_{X_{m_0+3Km}}\w))_{m\geq0}$ is a
sequence
of i.i.d.~random variables, under $P_0$, for any $m_0\geq K$. Therefore,
\begin{align*}
P_0\left(\forall m_0\ge K:\lim_{n\rightarrow\infty}
n^{-1}\!\!\!\!\!\!\!\sum_{3Km\leq n-m_0}
\!\!\!\!\!\!\!\!\Psi(T_{X_{m_0+3Km}}\w)=
\frac{E_0(\Psi(T_{X_{m_0}}\w))}{3K}\right)=1.
\end{align*}
It then follows that, $\P$-a.s.
\begin{align*}
P^\w_0\left(\lim_{n\rightarrow\infty}
n^{-1}\sum_{m=K}^{n-1}\Psi(T_{X_m}\w)=c\right)=1,
\end{align*}
where
$$c=(3K)^{-1}\sum_{m_0=K}^{4K-1}E_0(\Psi(T_{X_{m_0}}\w)).$$
Since the above quenched probability is ${\mathfrak S}_0$-measurable the convergence also
holds $\P_\infty$-a.s. But then $c$ cannot be anything other than $\E_\infty(\Psi)$.

By bounded convergence one then sees that
$n^{-1}\sum_{m=0}^{n-1}\Pi^m\Psi$ converges to $\E_\infty(\Psi)$,
$\P_\infty$-a.s. By $L^1$-approximations we get this same limit in
the $L^1$ sense for all $\Psi\in L^1(\P_\infty)$ and the
ergodicity follows from the development in Section IV.2 of
\cite{rosenblatt}. This proves Lemma \ref{ergodicity}.
\end{proof}

We continue with the proof of Theorem \ref{cltdirected}.
Next, we will show that condition (\ref{cond}) is satisfied.

\begin{lemma}
\label{bounddirected}
Under the assumptions of Theorem {\rm\ref{cltdirected}}
condition {\rm (\ref{cond})} is satisfied with $\alpha=1/4$.
\end{lemma}

\begin{proof}
Observe that
$$\norm{\sum_{k=0}^{n-1}\Pi^k g}^2_2=\sum_{i,j=0}^{n-1}
\sum_{x,y\in\Z^d}\int P_0^\w(X_i=x)P_0^\w(X_j=y)g(T_x\w)g(T_y\w)\P(d\w).$$
Note that the integral is taken with respect to $\P$ instead of $\P_\infty$.
This is because
all the integrands depend on $\w$ through $\mathfrak{S}_0$.

Now, since $\P$ is a product measure and $g$ has mean zero, the summands in the above
sum vanish unless
$i=j$ and $x=y$. 
Therefore, one has
$$\norm{\sum_{k=0}^{n-1}\Pi^k g}^2_2=\norm{g}^2_2\sum_{k=0}^{n-1}
\int P_{0,0}^\w(X_k=\tilde X_k)\P(d\w),$$
where $X_n$ and $\tilde X_n$ are two independent walkers driven by the same environment
$\w$. We will denote their law by $P_{0,0}^\w$.
To check condition (\ref{cond}) one needs to find the asymptotic
behaviour of the above sum.

Notice now that, under $\int P_{0,0}^\w \P(d\w)$, the difference $Y_n=X_n-\tilde X_n$ 
performs a random walk on $E$ with the following kernel:
\begin{align}
q(0,y)&=\sum_{z\in E}\E(\pi_{0,e_1+z}\pi_{0,e_1+z+y}),\nn\\
q(x,y)&=\sum_{z\in E}\E(\pi_{0,e_1+z})\E(\pi_{0,e_1+z+y-x}),~
\hbox{if }x\not=0.
\label{q}
\end{align}
This walk is actually a homogeneous symmetric random walk on $E$, perturbed
at $0$.  Due to Lemma 3.3 of \cite{FF} one then has
$$\sum_{k=0}^{n-1}\int P_{0,0}^\w(X_k=\tilde X_k)\P(d\w)=\Ord(\sqrt n).$$
Therefore, condition (\ref{cond}) is satisfied and we are done.
\end{proof}

\begin{remark}
Our application of Lemma 3.3 from \cite{FF} may appear unjustified
because of the additional hypotheses \cite{FF} employs. However,
(2.2) in \cite{FF} is superfluous. Once one notices that
$$q(0,y)>0\Rightarrow q_0(0,y)>0,$$
one can apply {\bf P7.1} on page 65 of Spitzer \cite{spitzer} to reduce the
treatment to a situation where (2.2) holds. Here, $q_0$ is the transition
kernel of the unperturbed random walk and is defined by (\ref{q}), for
all $x$. Also, in Lemma 3.2 of \cite{FF} the authors
reference {\bf P12.3} of Spitzer's book,
which requires more than just two moments on $q_0$. However, one can
instead use (3) of section 12 on page 122 of \cite{spitzer}. Furthermore,
Lemma 3.1 of \cite{FF} is not needed for our purposes, since
we only need the upper bound in (3.22) therein. Lastly, the reference
to {\bf P7.9} of \cite{spitzer}, in Lemma 3.3 of \cite{FF},
can be replaced by {\bf P7.6}. This allows to discard ``strong aperiodicity''.
\end{remark}

Now that we have verified all the assumptions of Theorem \ref{MW} are
satisfied, we can conclude the proof of Theorem \ref{cltdirected}.
Indeed, Theorem \ref{MW} implies that the claim of Theorem \ref{cltdirected}
holds for $\P_\infty$-a.e.~$\w$.
But since everything depends on $\w$ only through
$\mathfrak{S}_0$, and $\P_\infty=\P$ on $\mathfrak{S}_0$, the same holds for
$\P$-a.e.~$\w$.

Of course, Theorem \ref{MW} yields a different formula for the diffusion
matrix. However, the annealed invariance principle has to have
the same diffusion matrix which, as we have mentioned in the introduction,
is precisely $\mathfrak D$. We leave it for the reader to double-check, with a
direct calculation, that the two formulae do coincide.
\end{proof}

\begin{remark}
\label{BM}
Note that $(\bar X_1+M_1)\cdot e_1=0$ and, therefore, the Brownian motion in
question is actually $\nu$-dimensional or smaller. This is of course also
clear from the formula for $\mathfrak D$.
\end{remark}





\begin{appendix}
\section{An invariance principle for a vector-valued martingale
difference array.} 
\label{MG-CLT}
In this appendix we give a proof of the vector-valued martingale
invariance principle, needed in the proof of Theorem \ref{MW}, that is
based on the corresponding scalar result.
The scalar version appears as Theorem 7.4 in Chapter 7 of Durrett's textbook
\cite{durrett}. 
It is noteworthy that an invariance principle for general Banach space
valued martingale differences that unifies several results in
the literature has been proved by \cite{basu}. 

To avoid confusion with time $t$
we write in this section $A^T$ for the transpose of a vector or
matrix $A$, instead of the $A^t$ we have used in the rest of the
paper. An element of $\R^d$ is still regarded as a $d\times 1$
matrix.

Let  $(\Omega,\cG,P)$ be   a probability space on which are
defined sub-$\sigma$-algebras $\cG_{n,k}\subset\cG$ and
$\R^d$-valued random vectors $Y_{n,k}$. We say that
$\{Y_{n,k},\cG_{n,k}:n\ge 1,1\le k\le n\}$ is an $\R^d$-valued
square-integrable martingale difference array if the following properties
are satisfied:

\quad(i) $Y_{n,k}\hbox{ is }\cG_{n,k}\hbox{-measurable, }
\cG_{n,k-1}\subset\cG_{n,k},$

\quad(ii) $E(\abs{Y_{n,k}}^2)<\infty,$

\quad (iii) $E(Y_{n,k}\vert\cG_{n,k-1})=0,$\\ and in the last
condition we take $\cG_{n,0}=\{\phi,\Omega\}$. Define the
$\R^d$-valued processes $S_n(\cdot)$ by
\[S_n(t)=\sum_{k=1}^{[nt]} Y_{n,k}\] for $0\le t\le 1.$ The paths of
$S_n(\cdot)$ are in the Skorohod space $D_{\R^d}([0,1])$ of
$\R^d$-valued cadlag paths on $[0,1]$. Recall now from Remark
\ref{BMdef} the definition of a Brownian motion with diffusion
matrix $\Gamma$.  One then has the following:

\begin{theorem}
\label{mg-ip-thm} Let $\{Y_{n,k},\cG_{n,k}:n\ge 1,1\le k\le n\}$
be an $\R^d$-valued square-integrable martingale difference array
on a probability space $(\Omega,\cG,P)$. Let $\Gamma$ be a
symmetric, non-negative definite $d\times d$ matrix. Assume that
\begin{align}
\label{ass-1} \lim_{n\to\infty} \sum_{k=1}^{[nt]}
E(Y_{n,k}Y^T_{n,k}\vert\cG_{n,k-1})= t\Gamma \hbox{ in
probability,}
\end{align} for each $0\le t\le 1$, and
\begin{align}
\label{ass-2} \lim_{n\to\infty} \sum_{k=1}^{n}
E(\abs{Y_{n,k}}^2\ind\{\abs{Y_{n,k}}\ge \e\} \vert \cG_{n,k-1}) =
0 \hbox{ in probability,}
\end{align}
 for each $\e>0$.
Then $S_n(\cdot)$ converges weakly to a Brownian motion with
diffusion matrix $\Gamma$  on the Skorohod space
$D_{\R^d}([0,1])$.
\end{theorem}

\begin{proof}
The key to the proof is to apply a scalar martingale invariance
principle to one-dimensional projections of $S_n$, conditional on
the past.

Fix $0\le s<1$. Consider $k$  time points $0\le
s_1<s_2<\cdots<s_k\le s$, and a non-negative, bounded continuous
function $\Psi$ on $\R^{kd}$. Abbreviate
\[
Z_n=\Psi\left(S_n(s_1),S_n(s_2),\cdots,S_n(s_k)\right).
\]
Assume $E(Z_n)>0$ for all $n$. Pick also a non-zero vector
$\theta\in\R^d$, and a bounded continuous function $f$ on the
Skorohod space  $D_\R([0,1-s])$ of scalar-valued paths. Finally,
let $B_\theta$ denote a one-dimensional Brownian motion with
variance $E(B_\theta(t)^2)= \theta^T\Gamma\theta t$.

\begin{lemma}
\label{t-lm-1} We have the limit
\begin{align}
\label{lm-lim-1} \lim_{n\to\infty} \frac{E\left(f\left(\theta\cdot
S_n(s+\cdot)- \theta\cdot S_n(s)\right) Z_n\right)}{E(Z_n)}
=E(f(B_\theta)).
\end{align}
\end{lemma}

\begin{proof} 
Define a probability measure $\Ptil_n$ on $\Omega$ by
\[
\Ptil_n(A) = \frac1{E(Z_n)}E(\ind_A\cdot Z_n).
\]
$\Etil_n$ denotes expectation under $\Ptil_n$. Since $Z_n$ is
$\cG_{n,[ns]}$-measurable, we have
\[\Etil_n(h\vert\cG_{n,k})=E(h\vert\cG_{n,k})\] for any $k\ge [ns]$
and  $h\in L^1(P)$.

Define a scalar martingale difference array  $\{X_{n,m},\cF_{n,m}:
n\ge 1, 1\le m\le n-[ns]\}$ by
\[X_{n,m}= \theta\cdot Y_{n,[ns]+m}\quad\hbox{and}\quad
\cF_{n,m}=\cG_{n,[ns]+m}.\]
Observe first that by assumption (\ref{ass-1}) one has
\begin{align*}
V_{n,[nt]}&\equiv\sum_{j=1}^{[nt]}
\Etil_n(X_{n,j}^2\vert\cF_{n,j-1}) =\sum_{j=1}^{[nt]}
E(X_{n,j}^2\vert\cG_{n,[ns]+j-1})\\
&=\theta^T\left\{\sum_{k=1}^{[ns]+[nt]} E(Y_{n,k}Y^T_{n,k}\vert
\cG_{n,k-1})-\sum_{k=1}^{[ns]} E(Y_{n,k}Y^T_{n,k}\vert
\cG_{n,k-1}) \right\}\theta\\
&\longrightarrow t\theta^T\Gamma\theta \quad\hbox{in probability
as }n\to\infty.
\end{align*}
In case the reader is concerned that the first sum on the
second-last line above  goes up to $j=[ns]+[nt]$,
we point out that assumption (\ref{ass-2}) implies
\begin{align}
\label{max-0} \lim_{n\to\infty} \max_{1\le k\le n}\abs{Y_{n,k}}= 0
\quad\hbox{in probability.}
\end{align}
Thus the limits are not affected by finitely many terms.
(\ref{max-0}) follows from Dvoretsky's Lemma, by an argument that
can be found in the proof of Theorem (7.3) in Section 7.7 of
Durrett \cite{durrett} (see part (f) of that proof).

Next, by assumption (\ref{ass-2}), for any $\e>0$ there exists
$\e_0>0$ such that
\begin{align*}
\sum_{j=1}^{n-[ns]}
&\Etil_n(X_{n,j}^2\ind\{\abs{X_{n,j}}\ge\e\}\vert\cF_{n,j-1})\\
&\qquad\le\abs{\theta}^2 \sum_{k=1}^n E(\abs{Y_{n,k}}^2
\ind\{\abs{Y_{n,j}}\ge \e_0\}  \vert \cG_{n,k-1})\\
&\qquad\longrightarrow 0 \quad\hbox{in probability as
}n\to\infty.
\end{align*}

We have verified the hypotheses of the Lindeberg-Feller Theorem
for martingales that appears as Theorem (7.3) in Section 7.7 of
Durrett \cite{durrett}.  Consequently the process
\begin{align*}
U_n(t)&=\sum_{j=1}^{[nt]} X_{n,j}\\
&=\theta\cdot S_n(s+t)-\theta\cdot S_n(s)
-\theta\cdot Y_{n,[n(s+t)]} \ind\{[ns]+[nt]<[n(s+t)]\}
\end{align*}
satisfies
\begin{align}
\label{int-fU} \Etil_n(f(U_n))\to E(f(B_\theta)).
\end{align}
To be precise, Durrett's theorem treats  the continuous process
$\theta\cdot \Sbar_n(\cdot)$ defined by linear interpolation:
\[
\Sbar_n(t)=S_n(t) + (nt-[nt]) Y_{n,[nt]+1}\,,\quad 0\le t\le 1.
\]
But by (\ref{max-0}),
\begin{align}
\label{max-1} \sup_{0\le t\le 1} \abs{S_n(t)-\Sbar_n(t)} \to 0
\quad\mbox{in probability,}
\end{align}
so  the cadlag and continuous versions converge weakly together.

(\ref{int-fU})  is the same as (\ref{lm-lim-1}), again because by
(\ref{max-0}) whether $[ns]+[nt]$ differs from $[n(s+t)]$  is
immaterial for the limit. Lemma \ref{t-lm-1} is proved.
\end{proof}

Now we prove Theorem \ref{mg-ip-thm} from this lemma. First, by
taking $s=0$ and $\Psi\equiv 1$, $\theta\cdot\Sbar_n$  converges
weakly  to the Brownian motion $B_\theta$, for each vector
$\theta$. Thus
 all the scalar processes obtained as  projections of
$\{\Sbar_n\}$ are tight,  and hence   the vector-valued processes
$\{\Sbar_n\}$ themselves are tight on the space $C_{\R^d}([0,1])$.
And then by (\ref{max-1}),
  the vector-valued processes
$\{S_n\}$  are tight on the space $D_{\R^d}([0,1])$. This detour
via the continuous processes $\{\Sbar_n\}$ to get tightness of
$\{S_n\}$ was used because tightness of
 vector-valued  processes from projections is not
as obvious for cadlag paths  as it is for continuous paths (see
Exercise 22 from Chapter 3 of Ethier-Kurtz \cite{EK}).

Let a process $X$ be a weak limit point of $\{S_n\}$, and let
$S_{n_j}$ be the subsequence along which $S_{n_j}\Rightarrow X$.
The map $\eta\mapsto\theta\cdot\eta$ from $D_{\R^d}([0,1])$ into
$D_{\R}([0,1])$  is continuous, hence
 $\theta\cdot X$  has the distribution of the Brownian motion
$B_\theta$. It follows that $X$ has a version with almost surely
continuous paths. Then the finite-dimensional marginals converge
weakly:
\[
\left( S_{n_j}(s_1), S_{n_j}(s_2), \cdots, S_{n_j}(s_k)\right)
\Rightarrow \left( X(s_1), X(s_2), \cdots, X(s_k)\right).
\]
From all this we conclude that along the subsequence ${n_j}$ the
left-hand side of (\ref{lm-lim-1}) converges to
\[
\frac{E\left(f\left(\theta\cdot X(s+\cdot)- \theta\cdot
X(s)\right) \Psi\left( X(s_1),  \cdots, X(s_k)\right) \right)}
{E\left( \Psi\left( X(s_1),  \cdots, X(s_k)\right) \right)}.
\]
By  (\ref{lm-lim-1}) this must equal $E(f(B_\theta))$. Since the
time points $\{s_i\}$ and the function $\Psi$ are arbitrary, it
follows that
\begin{align}
\label{last-eq} E\left( e^{i\theta\cdot(X(s+t)-X(s))} \bigm\vert
X(r): 0\le r\le s\right)=E\left( e^{iB_\theta(t)}\right)
=e^{-\frac{t}2\theta^T\Gamma\theta}.~~
\end{align}
Varying the vector $\theta$ here implies that the increment
$X(s+t)-X(s)$ is independent  of the past up to time $s$, and is
distributed like the increment of a Brownian motion with diffusion
matrix $\Gamma$. Inductively on the number of increments we
conclude that $X$ has independent increments, continuous paths and
the correct Gaussian finite-dimensional distributions, which makes
it a Brownian motion with diffusion matrix $\Gamma$.

Note that it is critically important for this argument that in
(\ref{last-eq}) we can condition on $\{X(r): 0\le r\le s\}$ and
not only on
 $\{\theta\cdot X(r): 0\le r\le s\}$. This latter would
not suffice for the conclusion, as indicated by Exercise 2 in
Chapter 7 of \cite{EK}. This completes the proof of Theorem
\ref{mg-ip-thm}.
\end{proof}
\end{appendix}

\begin{acknowledgment}
The authors thank Dr.~M.~Bal\'azs and Prof.~S.R.S.~Varadhan for
valuable discussions.  Rassoul-Agha thanks the FIM for hospitality
and financial support during his visit to ETH, where part of this
work was performed. He also thanks Prof.~A-S.~Sznitman for
pointing out some useful references.
\end{acknowledgment}







\bibliographystyle{aop} 
\bibliography{refs}




%
%



\end{document}